\documentclass[a4paper,12pt]{article}
\usepackage[english]{babel}
\usepackage[T2A]{fontenc}
\usepackage[cp1251]{inputenc}
\usepackage[english]{babel}
\usepackage{amsthm}
\usepackage[tbtags]{amsmath}
\usepackage{amsfonts,amssymb}
\begin{document}

\newtheorem{theorem}{Theorem}
\newtheorem{lemma}{Lemma}
\newtheorem{proposition}{Proposition}
\newtheorem{Cor}{Corollary}

\begin{center}
{\bf Formal Matrix Rings: Isomorphism Problem}
\end{center}
\begin{center}
P.A. Krylov\footnote{National Research Tomsk State University,  e-mail: krylov@math.tsu.ru .},
A.A. Tuganbaev\footnote{National Research University <<MPEI>>, Lomonosov Moscow State University; e-mail: tuganbaev@gmail.com .}
\end{center}

\textbf{Abstract.} We consider the isomorphism problem for formal matrix rings over a given ring. Principal factor matrices of such rings play an important role in this case.
The work is supported by Russian Scientific Foundation, project 23-21-00375 (P.A. Krylov) and project 22-11-00052 (A.A. Tuganbaev).




\section{Introduction}\label{section1} 

Formal (or generalized) matrix rings over a given ring attract a lot of attention from specialists. It is natural, since such rings regularly appear in ring theory. In addition, they play an important role in the study of a number of classes of Artinian rings and algebras (see \cite{AusRS95}, \cite{BabO09}). They also serve as a source of varied examples for general ring theory. A number of aspects of the theory of formal matrix rings are presented in the book \cite{KryT17}.

There is one interesting type of formal matrix rings. In the case of $2\times 2$ matrices, they appeared in \cite{Kry08}. In the case of $n\times n$ matrices, they appeared in \cite{TanZ13}. We mean formal matrix rings over a given ring $R$ (or we say <<with values in $R$>>). This means that a concrete formal matrix ring has the same ring $R$ on all positions. The class of such rings is a direct expansion of an ordinary ring $M(n,R)$ of $n\times n$ matrices. However, properties of formal matrix rings over the ring $R$ may differ greatly from properties of the ring $M(n,R)$. In \cite[Chapter 4]{KryT17}, some questions are raised on formal matrix rings over a ring $R$ and three problems on formal matrix rings are formulated at the beginning of this chapter. In \cite[Sections 4.3--4.5]{KryT17}, these problems are solved for some types of formal matrix rings over the ring $R$. This paper is devoted to one of these three problems. Namely, this is isomorphism problem III. This problem is considered in \cite{AbyT15}, \cite{AbyT15b}, \cite{CheDS20}, \cite{Tap15}, \cite{Tap17}, \cite{Tap18}. 

In this paper, we consider only associative rings with non-zero identity element. If $R$ is a ring, then $M(n,R)$ is an ordinary ring of all $n\times n$ matrices with values in the ring $R$. The prime radical of an arbitrary ring $S$ is denoted by $P(S)$.

\section{Formal Matrix Rings over a Given Ring}\label{section2} 

We briefly recall the definition of a formal matrix ring. We fix a positive integer $n\ge 2$. Let $R_1,\ldots,R_n$ be rings and let $M_{ij}$ be $R_i$-$R_j$-bimodules with $M_{ii}=R_i$, $i,j=1,\ldots,n$. Let's assume that for any subscripts $i,j,k=1,\ldots,n$, we have an $R_i$-$R_k$-bimodule homomorphism $M_{ij}\otimes_{R_j}M_{jk}\to M_{ik}$. We denote by $K$ the set of all $n\times n$ matrices with values in bimodules $M_{ij}$. The set $K$ forms a ring with respect to standard matrix operations of addition and multiplication. Matrices are multiplies by the use of the above-mentioned bimodule homomorphisms. The ring $K$ is called a \textsf{formal} (or \textsf{generalized}) \textsf{matrix ring} of order $n$. The ring $K$ is of the following form:
$$
K=\begin{pmatrix}
R_1&M_{12}&\ldots&M_{1n}\\
M_{21}&R_{2}&\ldots&M_{2n}\\
\ldots&\ldots&\ldots&\ldots\\
M_{n1}&M_{n2}&\ldots&R_{n}
\end{pmatrix}.
$$
Let $R$ be some ring. If $K$ is a formal matrix ring such that $M_{ij}=R$ for all $i$ and $j$, then $K$ is called a \textsf{formal matrix ring over the ring} $R$ or a \textsf{formal matrix ring with values in the ring} $R$. Such rings can be defined directly. Namely, let
$\{s_{ijk}\,|\,i,j,k=1,\ldots,n\}$ be some set of central elements of the ring $R$ satisfying relations 
$$
s_{iik}=1=s_{ikk},\; s_{ijk}\cdot s_{ik\ell}=s_{ij\ell}\cdot s_{jk\ell} \eqno (1)
$$ 
for all subscripts $i,j,k,\ell=1,\ldots,n$. For arbitrary $n\times n$ matrices $A=(a_{ij})$ and $B=(b_{ij})$ with values in $R$, we define a new multiplication, by setting
$$
AB=C=(c_{ij}),\, \text{ where } c_{ij}=\sum_{k=1}^ns_{ikj}a_{ik}b_{kj}.
$$
As a result, we obtain a ring which is denoted by $K$ or $M(n,R,\Sigma)$, where $\Sigma$ is the set of all elements $s_{ijk}$. The set $\Sigma$ is called a \textsf{factor system} and its elements are called \textsf{factors} of the ring $K$. 
If all $s_{ijk}$ are equal to $1$, then we obtain the ordinary matrix ring $M(n,R)$.

The main relations $(1)$ imply the following relations:
$$
s_{iji}=s_{jij}=s_{ij\ell}\cdot s_{ji\ell}=s_{\ell ij}\cdot s_{\ell ji}.\eqno (2)
$$
It is useful to rewrite the last relation in $(2)$ in the form of three relations which follow from each other by permutation of subscripts:
\begin{gather*}
\begin{array}{c}
s_{iji}=s_{jij}=s_{ijk}\cdot s_{jik}=s_{kij}\cdot s_{kji},\\
s_{jkj}=s_{kjk}=s_{jki}\cdot s_{kji}=s_{ijk}\cdot s_{ikj},\\
s_{iki}=s_{kik}=s_{ikj}\cdot s_{kij}=s_{jik}\cdot s_{jki}.
\end{array}\qquad (3)
\end{gather*}

Let $\tau$ be a permutation of degree $n$. If $\Sigma=\{s_{ijk}\}$ is some factor system, then we set $t_{ijk}=s_{\tau(i)\tau(j)\tau(k)}$. Then $\{t_{ijk}\}$ is a factor system, as well, since it satisfies relations $(1)$. We denote it by $\tau\Sigma$. Consequently, there exists a formal matrix ring $M(n,R,\tau\Sigma)$. The rings $M(n,R,\Sigma)$ and $M(n,R,\tau\Sigma)$ are isomorphic to each other under the correspondence $A\to \tau A$, where $A=(a_{ij})$ and $\tau A=(a_{\tau(i)\tau(j)})$.

We can associate several matrices with a given ring $M(n,R,\Sigma)$. We set $S=(s_{iji})$ and $S_k=(s_{ikj})$ for every $k=1,\ldots,n$. These matrices are called \textsf{factor matrices} of the ring $M(n,R,\Sigma)$. The matrix $S$ is symmetrical. Following \cite{CheDS20}, we call it a \textsf{principal factor matrix}. In \cite{CheDS20} matrices $(s_{ijk})$ and $(s_{kij})$ are also used for $k=1,\ldots,n$. It is clear that the matrices $\tau S$ and $\tau S_k$ are the corresponding factor matrices for the ring $M(n,R,\tau\Sigma)$.

Sometimes, it is possible to choose a permutation $\tau$ such that the principal factor matrix $\tau S$ of the ring $M(n,R,\tau\Sigma)$ obtains a specific, simpler and more convenient form. We briefly present three corresponding cases. The first case slightly generalizes considerations in \cite[Section 4.3]{KryT17}; see Lemmas 4.3.1, 4.3.2 and the paragraph after the proof of Lemma 4.3.2.

Let $\Sigma$ be a factor system such that every factor $s_{ijk}$ is either non-invertible or equal to $1$. Under this assumption, additional relations appear between factors $s_{ijk}$. For example, with the use of relations $(3)$, it is easy to verify that the following lemma is true.

\textbf{Lemma 2.1.} Let $i,j,k$ be pairwise distinct subscripts.
Then elements $s_{iji}$, $s_{iki}$ and $s_{jkj}$ satisfy only one of the following three conditions.
\begin{enumerate}
\item[{\bf 1)}]
All three elements are equal to $1$.
\item[{\bf 2)}]
Some two elements of these three elements are non-invertible and the third element is equal to $1$.
\item[{\bf 3)}]
All three these elements are non-invertible.
\end{enumerate}

On the set $\{1,\ldots,n\}$, we define a binary relation $\sim$, by setting $i\sim j$ $\Leftrightarrow$ $s_{iji}=1$. Lemma 2.1 implies the following assertion.

\textbf{Lemma 2.2.} The relation $\sim$ is an equivalence relation.

Let's write the final result at the moment.

\textbf{Lemma 2.3.} There exists a permutation $\tau$ such that the matrix $\tau S$ can be presented in a block form such that blocks, consisting of $1$s, stay on the main diagonal and non-invertible elements stay on all remaining positions.

\textbf{Proof.} Let $\tau$ be a permutation such that in the upper row consists of numbers $1,\ldots,n$ in a natural order. The bottom row consists of equivalence classes of the relation $\sim$ which are arranged in random order. In every class, numbers are also arranged in random order. The matrix $\tau S$ has the structure specified in the lemma.~$\square$

Two other cases were considered later. One of these  cases is considered in \cite{KryN18} and \cite[Lemma 12.1, Lemma 12.2]{KryT21} and another case is considered in \cite[Lemma 4.1]{CheDS20}. Moreover, no restrictions are imposed on factors. However, the ring $R$ is assumed to be commutative in \cite{KryN18} and \cite{KryT21}. But in this context it doesn't matter, since factors are central elements.

Similarly to the above, omitting details, we can say that it turns out that \cite{KryN18} and \cite{KryT21} deal with the situation where always exists a  permutation $\tau$ such that the matrix $\tau S$ has the following block structure: the blocks on the main diagonal are filled with non-zero-divisors, and the remaining blocks are filled with zero-divisors.

In \cite[Lemma 4.1]{CheDS20}, it is proved that there exists a permutation $\tau$ such that the corresponding blocks on the main diagonal of the matrix $\tau S$ consist of invertible elements and all remaining blocks consist of non-invertible elements. The analogues of Lemma 2.1 and Lemma 2.2 are also true.

Of course, the three situations outlined can be combined within the framework of some general approach.

An interesting important class of formal matrix rings is formed by the rings $M(n,R,\Sigma)$ which factor systems consist of $1$s and some central element $s$. In \cite{CheDS20}, such systems $\Sigma$ are called \textsf{binary systems}. We denote the corresponding ring $M(n,R,\Sigma)$ by $M(n,R,s)$ and agree to call it the \textsf{$(s1)$-ring} of formal matrices.

The rings $M(n,R,s)$, where $s^2\ne 1$ and $s^2\ne s$, are studied in \cite[Sections 4.3, 4.4]{KryT17} and \cite{CheDS20}. 
Connections of such rings with crossed matrix rings are known; see \cite{BabO09}. In \cite[Lemma 4.7]{CheDS20}, it is proved that a principal factor matrix of the ring $M(n,R,s)$, where $s^2\ne 1$ and $s^2\ne s$, uniquely determines all remaining factor matrices. This result is very useful.

We obtain another quite specific type of rings if we set $s=0$. All factor matrices of such rings are $(01)$-matrices. Papers \cite{KryN18}, \cite{KryT21} and \cite{KryT22} contain various material on automorphism groups of formal matrix rings.

We return to an arbitrary ring $M(n,R,s)$. If the element $s$ is invertible, then we have an isomorphism $M(n,R,s)\cong M(n,R)$. Therefore, we can assume that the element $s$ is not invertible. Then it follows from Lemma 2.3 that there exists a permutation $\tau$ with the following property: the matrix $\tau S$ can be divided into blocks in such a way that the diagonal blocks are filled with $1$s, and all other blocks are filled with the element $s$. In such a situation, we say that the matrix $\tau S$ \textsf{has a canonical form}, and the matrix $S$ \textsf{is reduced to a canonical form}.

\textbf{Remark 2.4.} Of course, the canonical form is determined up to permutation of blocks on the main diagonal and the corresponding permutation of remaining blocks.

\section{Isomorphism Problem for Formal Matrix Rings}\label{section3}

\cite[Section 4.1]{KryT17} contains the isomorphism problem III for formal matrix rings:
\begin{itemize}
\item When do two factor systems define isomorphic formal matrix rings?
\end{itemize}
We consider this problem for $(s1)$-formal matrix rings.

We say that some ring $R$ \textsf{satisfies $(n,m)$-condition} if we have $m=n$ for any positive integers $n$ and $m$ such that the rings $M(n,R)$ and $M(m,R)$ are isomorphic to each other. For example, the $(n,m)$-condition holds if the ring $R$ is either commutative, or local, or is a left (right) principal ideal domain.

A ring $S$ is said to be \textsf{indecomposable} if $0$ and $1$ are only central idempotents of $S$.

\textbf{Theorem 3.1.} Let the factor ring $R/P(R)$ be indecomposable and satisfies $(n,m)$-condition and $s\in P(R)$. Let $K_1$ and $K_2$ be two $(s1)$-formal matrix rings with principal factor matrices $S$ and $T$, respectively. The following assertions hold.

\textbf{1.} If the rings $K_1$ and $K_2$ are isomorphic to each other, then the matrices $S$ and $T$ have the same canonical forms.

\textbf{2.} If $s^2\ne 1$ and $s^2\ne s$, then the converse is also true.

\textbf{Proof.} \textbf{1.} By Lemma 2.3, we can assume that the matrices $S$ and $T$ are presented in the canonical form. Let's assume that $K_1\cong K_2$. Then there exists a ring isomorphism
$$
\gamma\colon K_1/P(K_1)\to K_2/P(K_2).
$$
The structure of the prime radicals $P(K_1)$ and $P(K_2)$ is known (see \cite[Corollary 4.2.2]{KryT17} and the paragraph after the corollary). We also know the block structure of the matrices $S$ and $T$. With the use this information, we obtain relations
$$
K_1/P(K_1)=P_1\times\ldots\times P_k\, \text{ and }\, K_2/P(K_2)=Q_1\times\ldots\times Q_{\ell},
$$
where $k$ (resp., $\ell$) is the number of blocks on the main diagonal of the matrix $S$ (resp., $T$). In addition, all $P_i$ and $Q_j$ are full matrix rings of some orders. Since the ring $R/P(R)$ is indecomposable, all the rings $P_i/P(P_i)$ and $Q_j/P(Q_j)$ are indecomposable.

Here, we remark that there is an analogue of \cite[Lemma 9.6]{KryT21} (or  \cite[Lemma 6.1]{KryT22}) on automorphisms of direct products of indecomposable rings for isomorphisms between direct products of indecomposable rings. Therefore, $k=\ell$ and there exists a permutation $\tau$ of degree $k$ such that the restriction $\gamma$ to $P_i$ is an isomorphism $P_i\to Q_{\tau(i)}$, $i=1,\ldots,k$. Consequently, the canonical forms of the matrices $S$ and $T$ coincide.

\textbf{2.} Let $\{s_{ijk}\}$ (resp., $\{t_{ijk}\}$) be the set of all factors of the ring $K_1$ (resp., $K_2$). As noted earlier, factors of the form $s_{iji}$, i.e., elements of the the principal factor matrix $S$, determine all remaining factors $s_{ijk}$; the same is true for factors of the ring $K_2$. Thus, $s_{ijk}=t_{ijk}$ for all $i,j,k$. Therefore, we have the relation $K_1=K_2$.~$\square$

\textbf{Corollary 3.2.} The factor rings $K_1/P(K_1)$ and $K_2/P(K_2)$ are isomorphic to each other if and only if the matrices $S$ and $T$ have the same canonical form.

\textbf{Remark 3.3.} In \cite[Theorem 4.12]{CheDS20}, a result, which is similar to Theorem 3.1, is proved for a left Artinian ring $R$. 

The introduction to this article mentions isomorphism problem III from \cite[Section 4.1]{KryT17}. The following open question is a partial case of this problem.

\textbf{Open question.} Let $s$ and $t$ be two central elements of the ring $R$. When is the isomorphism $M(n,R,s)\cong M(n,R,t)$ true?

A similar question for some other rings $M(n,R,\Sigma)$ is considered in \cite[Section 4.5]{KryT17}.

\label{biblio}

\end{document}